\documentclass[amstex,12pt,russian,amssymb]{article}

\usepackage{mathtext}
\usepackage[cp1251]{inputenc}
\usepackage[T2A]{fontenc}
\usepackage[russian]{babel}
\usepackage[dvips]{graphicx}
\usepackage{amsmath}
\usepackage{amssymb}
\usepackage{amsxtra}
\usepackage{latexsym}
\usepackage{ifthen}

\begin{document}

\newcounter{lemma}
\newcommand{\lemma}{\par \refstepcounter{lemma}%
{\bf Лемма \arabic{lemma}.}}

\newcounter{corollary}
\newcommand{\corollary}{\par \refstepcounter{corollary}%
{\bf Следствие \arabic{corollary}.}}

\newcounter{remark}
\newcommand{\remark}{\par \refstepcounter{remark}%
{\bf Замечание \arabic{remark}.}}

\newcounter{theorem}
\newcommand{\theorem}{\par \refstepcounter{theorem}%
{\bf Теорема \arabic{theorem}.}}

\newcounter{proposition}
\newcommand{\proposition}{\par \refstepcounter{proposition}%
{\bf Предложение \arabic{proposition}.}}

\newcommand{\proof}{{\it Доказательство.\,\,}}
\renewcommand{\refname}{\centerline{\bf Список литературы}}

{\bf Р.Р.~Салимов} (Институт математики НАН Украины),

{\bf Е.А.~Севостьянов, А.А.~Маркиш} (Житомирский государственный
университет имени И.~Франко)

\medskip\medskip
\medskip
{\bf Р.Р.~Салімов} (Інститут математики НАН України),

{\bf Є.О.~Севостьянов, А.О.~Маркиш} (Житомирський державний
університет імені І.~Франко)

\medskip\medskip
\medskip
{\bf R.R.~Salimov} (Institute of Mathematics of NAS of Ukraine),

{\bf E.A.~Sevost'yanov, A.A.~Markysh} (Zhytomyr Ivan Franko State
University)

\medskip
{\bf Об оценке искажения расстояния снизу для одного класса
отображений }

\medskip
Изучается поведение одного подкласса отображений с конечным
искажением в окрестности начала координат. При определённых условиях
на характеристику квазиконформности установлена оценка искажения
расстояния снизу для таких отображений.

\medskip
{\bf Про оцінку спотворення відстані знизу для одного класу
відображень}

\medskip
Вивчається поведінка одного підкласу відображень зі скінченним
спотворенням у околі початку координат. За певних умов на
характеристику квазіконформності встановлено оцінку спотворення
відстані знизу для таких відображень.

\medskip
{\bf On lower estimate of distortion of a distance for one class of
mappings}

\medskip
A behavior of one class of mappings with finite distortion at a
neighborhood of the origin is investigated. There is proved a lower
estimate of distortion of a distance under mappings mentioned above.

\newpage
{\bf 1. Введение.} Настоящая заметка посвящена изучению отображений,
удовлетворяющих верхним оценкам искажения $p$-модуля семейства
кривых для $p\ne n.$ Напомним, что в работе \cite{SalSev$_1$} первых
двух авторов установлена оценка искажения расстояния при таких
отображениях, обобщающая классическую теорему К.~Икома
(см.~\cite{I}). Укажем также на публикации \cite{Sal$_1$} и
\cite{SalSev$_2$}, где, как и в \cite{SalSev$_1$}, присутствует
ограничение $n-1<p\leqslant n.$ На наш взгляд, случай $p>n$ также
заслуживает внимания, и именно он будет рассмотрен в данной заметке.
Как будет показано ниже, в случае $p>n$ мы имеем дело с оценкой
соответствующей величины снизу, а не сверху, что контрастирует с
ситуацией $n-1<p\leqslant n$. Указанное отличие связано с
принципиально иным поведением ёмкости при $p>n,$ на что указывает,
напр., неравенство~(8.8) в монографии~\cite{Ma}.

Стоит отметить, что классы исследуемых в работе отображений шире
традиционного <<конформного>> случая $p=n.$ Наиболее изученными
являются отображения с ограниченным искажением $p$-модуля. Ещё в
70-е годы минувшего столетия Ф.~Герингом установлена
квазиизометричность гомеоморфизмов, искажающих $p$-модуль в
ограниченное число раз при $n-1<p<n$ (см. \cite[теорема~3]{Ge}).
Поскольку квазиконформные отображения не обладают таким свойством,
данный факт указывает на целесообразность отдельного исследования
$p$-случая. Отметим также, что при $p>n$ гомеоморфизмы с аналогичным
свойством имеют квазиизометричные обратные отображения, в то время
как для прямых гомеоморфизмов этот факт, по-видимому, не установлен
(см. \cite[теорема~3]{Ge}).

\medskip
Напомним определения (см. \cite{MRSY}, \cite{GRSY} и \cite{GRY}).
Здесь и далее $D$ -- область в ${\Bbb R}^n,$ $n\geqslant 2,$ $m$ --
мера Лебега в ${\Bbb R}^n,$ отображение $f:D\rightarrow {\Bbb R}^n,$
$x=(x_1,\ldots, x_n),$ $f(x)=(f_1(x),\ldots, f_n(x)),$
предполагается непрерывным. Напомним, под что семейством кривых
$\Gamma$ подразумевается некоторый фиксированный набор кривых
$\gamma,$ а $f(\Gamma)=\left\{f\circ\gamma|\gamma\in\Gamma\right\}.$
Всюду далее
$$A(x_0,r_1,r_2)=\{ x\,\in\,{\Bbb R}^n :
r_1<|x-x_0|<r_2\}\,, S(x_0,r)=\{x\in {\Bbb R}^n: |x-x_0|=r\}\,,$$
$$B(x_0, r)=\left\{x\in{\Bbb R}^n: |x-x_0|< r\right\},\quad {\Bbb
B}^n := B(0, 1)\,, \quad {\Bbb S}^{n-1}:=S(0, 1)\,,$$
$\Omega_n$ -- объём единичного шара ${\Bbb B}^n$ в ${\Bbb R}^n,$ а
$\omega_{n-1}$ --  площадь единичной сферы ${\Bbb S}^{n-1}$ в ${\Bbb
R}^n.$ Для произвольных множеств $E,$ $F\subset \overline{{\Bbb
R}^n}:={\Bbb R}^n\cup\{\infty\}$ через $\Gamma(E,F,D)$ в дальнейшем
обозначается семейство всевозможных кривых
$\gamma:[a,b]\rightarrow\overline{{\Bbb R}^n},$ соединяющих $E$ и
$F$ в $D$ (т.е., $\gamma(a)\in E,$ $\gamma(b)\in F$ и $\gamma(t)\in
D$ при $t\in (a, b)$). Следующие определения могут быть найдены,
напр., в \cite[разд.~1--6, гл.~I]{Va}. Борелева функция $\rho:{\Bbb
R}^n\,\rightarrow [0,\infty]$ называется {\it допустимой} для
семейства $\Gamma$ кривых $\gamma$ в ${\Bbb R}^n,$ если
$\int\limits_{\gamma}\rho (x)\,|dx|\geqslant 1$
%
%\end{equation}
%
для всех кривых $ \gamma \in \Gamma.$ В этом случае мы пишем: $\rho
\in {\rm adm} \,\Gamma .$
Пусть $p\geqslant 1,$ тогда {\it $p$ -- модулем} семейства кривых
$\Gamma $ называется величина
$M_p(\Gamma)=\inf\limits_{\rho \in \,{\rm adm}\,\Gamma}
\int\limits_D \rho ^p (x)\ \ dm(x)\,.$
Пусть $x_0\in D,$ $r_0\,=\,{\rm dist}\, (x_0,\partial D),$
$Q:D\rightarrow\,[0,\infty]$ -- некоторая заданная измеримая по
Лебегу функция. Обозначим через $S_i:=S(x_0, r_i),$ где
$0<r_1<r_2<\infty.$ Предположим, что отображение $f$ удовлетворяет
для каждых $0<r_1<r_2< r_0={\rm dist}\,(x_0,
\partial D)$ условию
$$
M_p\left(f\left(\Gamma\left(S_1,\,S_2,\,A\right)\right)\right)\leqslant
\int\limits_{A} Q(x)\cdot \eta^p(|x-x_0|) \,dm(x)\,,$$
выполненному для произвольной измеримой по Лебегу функции $\eta:
(r_1,r_2)\rightarrow [0,\infty ]$ такой, что
$ \int\limits_{r_1}^{r_2}\eta(r) dr\geqslant 1,$
где $A=A(x_0,r_1,r_2)$ -- сферическое кольцо с центром в точке $x_0$
радиусов $r_1$ и $r_2.$ Тогда будем говорить, что $f$ является {\it
кольцевым отображением в точке $x_0\,\in\,D$ относительно
$p$-модуля.} В настоящей работе нами устанавливается справедливость
следующего результата.

\medskip
\begin{theorem}\label{th1}{\sl\,
Пусть $n\geqslant 2,$ $n<p<\infty,$ $f:{\Bbb B}^n\rightarrow {\Bbb
R}^n$ -- открытое дискретное кольцевое $Q$-отоб\-ра\-же\-ние
относительно $p$-модуля в точке $x_0=0,$ $f(0)=0.$ Предположим, что
$Q:{\Bbb B}^n\rightarrow[0,\,\infty]$ -- локально интегрируемая
функция в ${\Bbb B}^n,$ удовлетворяющая условию
\begin{equation}\label{eq1A}
Q_0=\liminf\limits_{\varepsilon\rightarrow
0}\frac{1}{\Omega_n\cdot\varepsilon^n} \int\limits_{B(0,
\varepsilon)} Q(x)\,dm(x)>0\,.
\end{equation}
Тогда имеет место оценка:
$$
\limsup\limits_{x\rightarrow 0}\frac{|f(x)|}{|x|}\geqslant c_{0}
Q_0^{\frac{1}{n-p}}\,,
$$
где  $c_0$  -- некоторая положительная постоянная, зависящая только
от размерности пространства $n$ и $p.$}
\end{theorem}

\medskip
{\bf 2. Вспомогательные результаты. Оценка верхнего предела одной
функции. } Пару $E=(A,C),$ где $A\subset{\Bbb R}^n$ -- открытое
множество и $C$ -- непустое компактное множество, содержащееся в
$A,$ будем называть {\it конденсатором} в ${\Bbb R}^n.$ Говорят
также, что конденсатор $E=(A,C)$ лежит в области $D,$ если $A\subset
D.$ Очевидно, что если $f:D\rightarrow{\Bbb R}^n$ -- непрерывное
открытое отображение и $E=(A,C)$ -- конденсатор в $D,$ то
$f(E):=(f(A),f(C))$ также является конденсатором в $f(D).$

\medskip
Обозначим через $C_0(A)$  множество всех непрерывных функций
$u:A\rightarrow{\Bbb R}^1$ с компактным носителем, $W_0(E)=W_0(A,C)$
-- семейство неотрицательных функций $u:A\rightarrow{\Bbb R}^1$
таких, что 1) $u\in C_0(A),$ 2) $u(x)\geqslant1$ для $x\in C$ и 3)
$u$ принадлежит классу $ACL.$ Также обозначим
$ |\nabla u|={\left(\sum\limits_{i=1}^n\,{\left(\frac{\partial
u}{\partial x_i}\right)}^2 \right)}^{1/2}.$
При $p\geqslant 1$ величину
$${\rm cap}_p\,E={\rm cap}_p\,(A,C)=\inf\limits_{u\in W_0(E)}\,
\int\limits_{A}\,|\nabla u|^p\,dm(x)$$
называют {\it $p$-ёмкостью} конденсатора $E.$ При $n<p<\infty$
\begin{equation}\label{maz} {\rm cap}_p\,E={\rm cap}_p\,(A, C)\geqslant
n{\Omega_n}^{\frac{p}{n}}
\left(\frac{p-n}{p-1}\right)^{p-1}\left((m(A))^{\frac{p-n}{n(p-1)}}-
(m(C))^{\frac{p-n}{n(p-1)}}\right)^{1-p}\,,\end{equation}
где ${\Omega}_n$ -  объем  единичного шара  в ${\Bbb R}^n,$ см.,
напр., \cite[неравенство~(8.7)]{Ma}.

\medskip
Пусть $Q:D\rightarrow[0,\infty]$ -- измеримая по Лебегу функция.
Тогда положим
$$q_{x_0}(r)=
\frac{1}{\omega_{n-1}r^{n-1}}\int\limits_{S(x_0,\,r)}Q(x)d\mathcal{H}^{n-1}\,,$$
где $\mathcal{H}^{n-1}$ -- $(n-1)$-мерная мера Хаусдорфа. Следующая
лемма при $p\in (1, n]$ доказана в \cite[лемма~1]{SalSev$_1$}. В
случае произвольного $p>1$ её доказательство дословно повторяет
рассуждения, относящиеся к случаю $p\in (1, n],$ и потому
опускается.

\medskip
\begin{lemma}\label{lem1} {\sl Пусть $n\geqslant 2,$ $p\geqslant 1,$ $Q:D\rightarrow [0, \infty]$ --
заданная измеримая по Лебегу функция, $f:D\rightarrow
\overline{{\Bbb R}^n}$ -- открытое дискретное кольцевое
$Q$-отображение в точке $x_0\in D$ относительно $p$-модуля и $E$ --
конденсатор вида $E=\left(B(x_0, r_2), \overline{B(x_0,
r_1)}\right),$ $0<r_1<r_2< {\rm dist} \, (x_0,\partial D).$ Полагаем
\begin{equation}\label{eq9}
I\ =\ I(x_0,r_1,r_2)\ =\ \int\limits_{r_1}^{r_2}\
\frac{dr}{r^{\frac{n-1}{p-1}}q_{x_0}^{\frac{1}{p-1}}(r)}\,.
\end{equation}
Тогда для конденсатора $f(E)=\left(f\left(B(x_0, r_2)\right),
f\left(\overline{B(x_0, r_1)}\right)\right)$ выполнено соотношение
%
%\begin{equation}\label{eq2A}
$${\rm cap}_p\, f(E)\leqslant\frac{\omega_{n-1}}{I^{p-1}}\,.$$
%
%\end{equation}
}
\end{lemma}

\medskip
Аналог следующей леммы доказан в \cite[лемма~5]{SalSev$_1$}.

\medskip
\begin{lemma}\label{prA3} {\sl\, Пусть $f:{\Bbb B}^n \rightarrow {\Bbb R}^n,$ $n\geqslant 2,$
-- открытое отображение, удовлетворяющее условию $f(0)=0.$
Предположим, что существует функция $R: [0, 1]\rightarrow [0,
\infty),$ такая что
\begin{equation}\label{eqA25} m\left(f(B(0, r))\right)\geqslant \Omega_n\,R^{n}(r)\,.\end{equation}
Тогда
$$\limsup\limits_{x\rightarrow 0}\frac{|f(x)|}{R(|x|)}\geqslant
1\,.$$}
\end{lemma}

\begin{proof} Полагаем $\max\limits_{|x|=r}|f(x)|=L_f(r).$
Покажем, что \begin{equation}\label{eqA27} f\left(B(0,
r)\right)\subset B(0, L_f(r))
\end{equation}
при каждом $r\in (0, 1).$ Для этого зафиксируем $r_0\in (0, 1)$ и
обозначим $M:=\sup\limits_{y\in \overline{f(B(0, r_0))}}|y|.$ По
определению точной верхней грани, найдётся последовательность
$y_k\in f(B(0, r_0)),$ $|y_k|\rightarrow M$ при
$k\rightarrow\infty.$ Тогда $y_k=f(x_k),$ $x_k\in B(0, r_0).$ Так
как $\overline{B(0, r_0)}$ -- компакт в ${\Bbb B}^n,$ мы можем
считать, что при некотором $x_0\in \overline{B(0, r_0)}$ выполнено
условие $x_k\rightarrow x_0$ при $k\rightarrow\infty.$ Поскольку $f$
-- непрерывное отображение в ${\Bbb B}^n,$ то $f(x_k)\rightarrow
f(x_0)$ при $k\rightarrow\infty,$ так что $f(x_0)=y_0.$ Таким
образом, $y_0\in f(\overline{B(0, r_0)}).$ Значит,
\begin{equation}\label{eq1}
M:=\max\limits_{y\in \overline{f(B(0, r_0))}}|y|=|y_0|\,,\quad
y_0\in f(\overline{B(0, r_0)})\,.
\end{equation}
Заметим, что в силу открытости отображения $f$ случай $y_0\in f(B(0,
r_0))$ невозможен. В самом деле, если бы $y_0\in f(B(0, r_0)),$ то
тогда $y_0$ входило бы во множество $f(B(0, r_0))$ вместе с
некоторой своей окрестностью $B(y_0, \delta),$ кроме того, $y_0\ne
0$ ввиду открытости отображения $f.$  Представим $y_0$ в виде:
$y_0=|y_0|\cdot\frac{y_0}{|y_0|}.$ Тогда вектор
$\widetilde{y_0}:=(|y_0|+\delta/2)\cdot\frac{y_0}{|y_0|}$ имеет
модуль больший, чем $y_0$ и всё ещё лежит в $f(B(0, r_0)).$ Однако,
последнее противоречит определению $y_0.$ Полученное противоречие
указывает на то, что $y_0\not\in f(B(0, r_0))$ и, значит, $y_0\in
\partial f(B(0, r_0)).$ В частности, отсюда следует, что
\begin{equation}\label{eq2}
|f(x)|<M\qquad \forall\,x\in B(0, r_0)\,.
\end{equation}
Поскольку $y_0\in
\partial f(B(0, r_0)),$  в силу открытости отображения $f$ имеем: $y_0\in
f(S(0, r_0)).$

\medskip
Итак, $y_0=f(x_0),$ где $x_0\in S(0, r_0).$ В таком случае, ввиду
соотношений (\ref{eq1}) и (\ref{eq2}) для всякого $x\in B(0, r_0)$
мы получим, что
$$|f(x)|< M=|y_0|=|f(x_0)|\leqslant \sup\limits_{x\in S(0, r_0)}|f(x)|=L_f(r_0)\,,$$
так что $f(x)\in B(0, L_f(r_0)).$ Включение (\ref{eqA27})
установлено.

Из соотношения (\ref{eqA27}), учитывая условие $f(0)=0,$ имеем
$\Omega_{n}\,L_{f}^{n}(r)\geqslant m(f(B(0, r)))$ и, следовательно,
\begin{equation}\label{eqA28}
L_{f}(r)\geqslant \left(\frac{m(f(B(0,
r)))}{\Omega_n}\right)^{\frac{1}{n}}\,.\end{equation}
Таким образом, учитывая  неравенства (\ref{eqA25}) и (\ref{eqA28}),
получаем
$$\limsup\limits_{x\rightarrow 0}\frac{|f(x)|}{R(|x|)}=\limsup\limits_{r\rightarrow 0}\frac{L_{f}(r)}{R(r)}
\geqslant \limsup\limits_{r\rightarrow 0}\left(\frac{m(f(B(0,
r)))}{\Omega_n}\right)^{\frac{1}{n}}\cdot\frac{1}{R(r)}\geqslant
1\,.$$
Лемма \ref{prA3} доказана. \end{proof}$\Box$

\medskip
{\bf 3. Доказательство основного результата.} {\it Доказательство
теоремы \ref{th1}.} Рассмотрим кольцо $A=A(0, \varepsilon_1,
\varepsilon_2),$ $0<\varepsilon_1<\varepsilon_2<1.$ Пусть $E$ --
конденсатор вида $E=\left(B(0, \varepsilon_2), \overline{B(0,
\varepsilon_1)}\right).$ Положим
$$
 \eta_0(r)=\frac{1}{Ir^{\frac{n-1}{p-1}}q_{x_0}^{\frac{1}{p-1}}(r)}\,,
$$
где $I$ -- величина, определённая в (\ref{eq9}). Согласно
\cite[лемма~2.2]{Sal})
\begin{equation}\label{eq10}
\frac{\omega_{n-1}}{I^{p-1}}\ =\ \int\limits_{A} Q(x)\cdot
\eta_0^p(|x-x_0|)\ dm(x)\ \leqslant\ \int\limits_{A} Q(x)\cdot
\eta^p(|x-x_0|)\ dm(x)
\end{equation}
для фиксированной измеримой функции $Q:{\Bbb R}^n\,\rightarrow
[0,\infty]$ такой, что $q_{x_0}(r)\ne \infty$ для п.в. $r>0,$ и
любой функции $\eta :(r_1,r_2)\rightarrow [0,\infty]$ такой, что
$ \int\limits_{r_1}^{r_2}\eta(r) dr = 1.$
Ввиду леммы \ref{lem1} и соотношения (\ref{eq10}) неравенство
\begin{equation}\label{eq12}
{\rm cap}_p\,f(E)\leqslant \int\limits_{A} Q(x)\cdot
\eta^p(|x-x_0|)\, dm(x)
\end{equation}
будет выполнено для произвольной измеримой по Лебегу функции $\eta:
(\varepsilon_1, \varepsilon_2)\rightarrow [0,\infty ]$ такой, что
$\int\limits_{\varepsilon_1}^{\varepsilon_2}\eta(r) dr\geqslant 1.$
Заметим, что функция
$$ \eta(t)\,=\,\left
\{\begin{array}{rr} \frac{1}{\varepsilon_2-\varepsilon_1}, &  \ t\in (\varepsilon_1,\varepsilon_2) \\
0, & \ t\in \Bbb{R}\setminus (\varepsilon_1,\varepsilon_2)
\end{array}\right.
$$
удовлетворяет условию
$\int\limits_{\varepsilon_1}^{\varepsilon_2}\eta(r) dr\geqslant 1,$
поэтому, согласно (\ref{eq12}) мы получим, что
$${\rm cap}_p\,\left(f(B(0, \varepsilon_2)),f(\overline{B(0,
\varepsilon_1)})\right) \leqslant
\frac{1}{(\varepsilon_2-\varepsilon_1)^p}\int\limits_{A(0,
\varepsilon_1, \varepsilon_2)} Q(x)\ dm(x)\ .$$
Далее, выбирая $\varepsilon_1=\varepsilon$ и
$\varepsilon_2=2\varepsilon$, получим
\begin{equation}\label{eq101}{\rm cap}_p\,\left(f(B(0, 2\varepsilon)),f(\overline{B(0, \varepsilon)})\right)\leqslant\,
\frac{1}{\varepsilon^p}\int\limits_{A(0, \varepsilon,
2\varepsilon)}Q(x)\,dm(x)\,.
\end{equation}
С другой стороны,  в силу  неравенства (\ref{maz}) при каждом
фиксированном $\varepsilon>r>0$ мы имеем оценку:
$${\rm cap}_p\,\left(f(B(0, 2\varepsilon)),f(\overline{B(0,
\varepsilon)})\right)\geqslant {\rm cap}_p\,\left(f(B(0,
2\varepsilon)),f(\overline{B(0, r)})\right)\geqslant$$
\begin{equation}\label{eq102}
\geqslant n{\Omega_n}^{\frac{p}{n}}
\left(\frac{p-n}{p-1}\right)^{p-1}\left((m(f(B(0,
2\varepsilon))))^{\frac{p-n}{n(p-1)}}- (m(f(\overline{B(0,
r)})))^{\frac{p-n}{n(p-1)}}\right)^{1-p}\,.
\end{equation}
Соотношение (\ref{eq102}) имеет место при любом $r\in (0,
\varepsilon),$ поэтому можно перейти к пределу при $r\rightarrow 0.$
В таком случае, мы получим:
\begin{equation}\label{eq3}
{\rm cap}_p\,\left(f(B(0, 2\varepsilon)),f(\overline{B(0,
\varepsilon)})\right)\geqslant c\cdot (m(f(B(0,
2\varepsilon))))^{\frac{n-p}{n}}\,,
\end{equation}
где $c:=n^{1/(1-p)}{\Omega_n}^{p/(n(1-p))} \cdot (p-1)/(p-n).$
Комбинируя   (\ref{eq101}) и (\ref{eq3}), мы получаем, что
\begin{equation}\label{eq4.2} \frac{m(f(B(0, 2\varepsilon))}
{2^n\Omega_n \varepsilon^n }\geqslant c_1 \
\left(\frac{1}{2^n\Omega_n\cdot\varepsilon^n}\int\limits_{B(0,
2\varepsilon)} Q(x)\, dm(x)\right)^{\frac{n}{n-p}} \,,\end{equation}
где $c_1$ - положительная постоянная зависящая только от $n$ и $p.$

Положим $L_f(\varepsilon)=\max\limits_{|x|=\varepsilon}|f(x)|.$
Используя соотношение (\ref{eqA28}), мы получим, что
\begin{equation}\label{eqks1.200}
\limsup\limits_{x\rightarrow
0}\frac{|f(x)|}{|x|}=\limsup\limits_{\varepsilon\rightarrow
0}\frac{L_f(2\varepsilon)}{2\varepsilon}
\geqslant\limsup\limits_{\varepsilon\rightarrow
0}\left(\frac{m(f(B(0,
2\varepsilon))}{\Omega_n(2\varepsilon)^n}\right)^{\frac{1}{n}}\,.
\end{equation}
Наконец, комбинируя (\ref{eq4.2}) и (\ref{eqks1.200}), имеем:
$$
\limsup\limits_{x\rightarrow 0}\frac{|f(x)|}{|x|}\geqslant c_{0} \
\limsup\limits_{\varepsilon\rightarrow
0}\left(\frac{1}{2^n\Omega_n\cdot\varepsilon^n}\int\limits_{B(0,
2\varepsilon)} Q(x)\, dm(x)\right)^{\frac{1}{n-p}}=c_0\cdot
Q_0^{\frac{1}{n-p}}\,,
$$
где $c_0>0$ -- некоторая постоянная, зависящая только от $n$ и $p.$
Теорема  доказана. $\Box$

\medskip
Из теоремы~\ref{th1} получаем следующее (даже более общее)
утверждение.

\medskip
\begin{corollary}\label{cor2}{\sl\,
Пусть $D$ -- область в ${\Bbb R}^n,$ $n\geqslant 2,$ $n<p<\infty,$
$f:D\rightarrow {\Bbb R}^n$ -- открытое дискретное кольцевое
$Q$-отоб\-ра\-же\-ние относительно $p$-модуля в точке $x_0,$ $x_0\in
D.$ Предположим, что $Q:D\rightarrow[0,\,\infty]$ -- локально
интегрируемая функция, удовлетворяющая условию
\begin{equation}\label{eq1B}
Q_0=\liminf\limits_{\varepsilon\rightarrow
0}\frac{1}{\Omega_n\cdot\varepsilon^n} \int\limits_{B(x_0,
\varepsilon)} Q(x)\,dm(x)>0\,.
\end{equation}
Тогда имеет место оценка:
\begin{equation}\label{eq15}
\limsup\limits_{x\rightarrow
x_0}\frac{|f(x)-f(x_0)|}{|x-x_0|}\geqslant c_0
Q_0^{\frac{1}{n-p}}\,,
\end{equation}
где  $c_0$  -- некоторая положительная постоянная, зависящая от
точки $x_0,$ размерности пространства $n$ и числа $p.$}
\end{corollary}

\medskip
{\it Доказательство} следствия~\ref{cor2} легко вытекает из
теоремы~\ref{th1}. В самом деле, если имеется кольцевое
$Q$-отображение $f:D\rightarrow {\Bbb R}^n$ в точке $x_0\in D,$ то
зафиксируем $r>0$ так, чтобы $\overline{B(x_0, r)}\subset D.$
Заметим, что вспомогательное преобразование
$\widetilde{f}(y):=f(ry+x_0)-f(x_0)$ также является
$\widetilde{Q}$-кольцевым отображением $f:{\Bbb B}^n\rightarrow{\Bbb
R}^n$ относительно $p$-модуля в точке 0, где
$\widetilde{Q}(y)=r^{n-p}\cdot Q(ry+x_0)$ и, кроме того, $f(0)=0$
(см., напр., \cite[теорема~8.2]{Va} по поводу изменения $p$-модуля
при растяжениях). Кроме того, заметим, что условие (\ref{eq1B})
выполняется для функции $\widetilde{Q}$ при $\widetilde{Q}_0$ вместо
$Q_0,$ где $\widetilde{Q}_0=r^{n-p}\cdot Q_0.$ Из этой теоремы
следует, что
$$
\limsup\limits_{y\rightarrow 0}\frac{|\widetilde{f}(y)|}{|y|}=
\limsup\limits_{y\rightarrow 0}\frac{|
f(ry+x_0)-f(x_0)|}{|y|}\geqslant c_0 Q_0^{\frac{1}{n-p}}\cdot r\,.
$$
Если в последнем соотношении перейти к переменной $x:=ry+x_0,$ то
оно запишется в виде~(\ref{eq15}), что и требовалось
установить.~$\Box$

\medskip
В качестве ещё одного полезного следствия из теоремы~\ref{th1} имеем
следующее утверждение.

\medskip
\begin{corollary}\label{cor1}
{\sl Предположим, что в условиях теоремы~\ref{th1} выполнено
равенство $Q_0=0.$ Тогда $\limsup\limits_{x\rightarrow
0}\frac{|f(x)|}{|x|}=\infty.$}
\end{corollary}

\medskip
\begin{proof}
Повторяя рассуждения, приведенные при доказательстве
леммы~\ref{lem1}, на основании соотношения вида (\ref{eqA28}), будем
иметь:
\begin{equation}\label{eq13}
L_f(r)\geqslant \left(\frac{m(f(B(0,
r)))}{\Omega_n}\right)^{\frac{1}{n}}\,,\end{equation}
где, как и ранее, $L_f(r):=\sup\limits_{x\in S(0, r)}|f(x)|.$
Повторяя также рассуждения, приведённые при доказательстве
теоремы~\ref{th1}, мы получим соотношение вида (\ref{eq4.2}), из
которого следует, что
\begin{equation}\label{eq14}
\lim\limits_{\varepsilon\rightarrow 0} \frac{m(f(B(0,
2\varepsilon))} {2^n\Omega_n \varepsilon^n }=\infty\,.
\end{equation}
Тогда из (\ref{eq13}) и (\ref{eq14}) получаем:
$$\limsup\limits_{x\rightarrow 0}\frac{|f(x)|}{|x|}=\limsup\limits_{r\rightarrow 0}
\frac{L_f(r)}{r}=\infty\,,$$
что и требовалось установить.~$\Box$
\end{proof}

\medskip
{\bf Пример.} Следующий пример указывает на содержательность условий
и заключения следствия~\ref{cor1}. При фиксированном $\theta\in (0,
1)$ рассмотрим отображение
$$f(x)=\frac{x}{|x|}|x|^{\theta},\quad x\in {\Bbb B}^n\,.$$
Заметим, что $|f(x)|/|x|\rightarrow\infty$ при $x\rightarrow 0.$
Покажем, что $f$ -- кольцевое $Q$-отображение в нуле. Для этого
воспользуемся \cite[теорема~2.2]{SalSev$_3$}. Очевидно, $f\in
C^1({\Bbb B}^n\setminus\{0\}),$ так что $f$ дифференцируемо почти
всюду в ${\Bbb B}^n,$ принадлежит классу $W^{1, p}_{loc}({\Bbb
B}^n\setminus\{0\})$ и обладает $N$-свойством. По тем же причинам
отображение $f^{\,-1}(y)=\frac{|y|^{1/\theta}}{|y|}y$ также
принадлежит классу $W^{1, p}_{loc}({\Bbb B}^n\setminus\{0\})$ и
обладает $N$-свойством. По \cite[теорема~2.2]{SalSev$_3$}
отображение $f$ является кольцевым $Q$-отображением в каждой точке
$x_0\in{\Bbb B}^n$ при
$$Q(x):=K_{I, p}(x,f)\quad =\quad\left\{
\begin{array}{rr}
\frac{J(x,f)}{l^p(x, f)}, & J(x,f)\ne 0,\\
1,  &  f^{\,\prime}(x)=0, \\
\infty, & \text{в\,\,остальных\,\,случаях},
\end{array}
\right.$$
где $l\left(f^{\,\prime}(x)\right)\,=\,\,\,\min\limits_{h\in {\Bbb
R}^n \setminus \{0\}} \frac {|f^{\,\prime}(x)h|}{|h|}$ и $J(x,
f)={\rm det}\, f^{\,\prime}(x).$ Подсчитаем $K_{I, p}(x,f),$ для
чего воспользуемся предложением 5.1 в \cite{IS}. В обозначениях
этого предложения
$$\lambda_r(x)=\frac{\partial |f(x)|}{\partial|x|}=\theta|x|^{\theta-1}\,,\quad\lambda_{\tau}(x)=\frac{|f(x)|}{|x|}=|x|^{\theta-1}\,.$$
Заметим, что $l\left(f^{\,\prime}(x)\right)=\min\{\lambda_r(x),
\lambda_{\tau}(x)\}$ и $|J(x,
f)|=\lambda^{n-1}_{\tau}(x)\cdot\lambda_r(x)$ (см.
\cite[разд.~5.1]{IS}). Очевидно, $\lambda_{r}(x)<\lambda_{\tau}(x)$
и
$$Q(x):=K_{I, p}(x, f)=\frac{\theta|x|^{(n-1)(\theta-1)}\cdot |x|^{\theta-1}}{(\theta|x|^{\theta-1})^p}=
\theta^{1-p}\cdot|x|^{(n-p)(\theta-1)}\,.$$
Вычислим $Q_0$ по формуле (\ref{eq1A}). По теореме Фубини будем
иметь:
$$\int\limits_{B(0, \varepsilon)}Q(x)\,dm(x)=\int\limits_0^{\varepsilon}\int\limits_{S(0, r)}
\theta^{1-p}\cdot|x|^{(n-p)(\theta-1)}\,dS
dr=\omega_{n-1}\theta^{1-p}\int\limits_{0}^{\varepsilon}r^{n-1}\cdot
r^{(n-p)(\theta-1)}dr=$$
$$=\frac{\omega_{n-1}\theta^{1-p}\varepsilon^{(n-p)(\theta-1)+n}}{(n-p)(\theta-1)+n}=C\cdot \varepsilon^{(n-p)(\theta-1)+n}\,,$$
где $C:=\frac{\omega_{n-1}\theta^{1-p}}{(n-p)(\theta-1)+n}.$
Учитывая полученное выше, имеем:
$$Q_0:=\limsup\limits_{\varepsilon\rightarrow 0}\frac{C}{\Omega_n}\cdot
\frac{\varepsilon^{(n-p)(\theta-1)+n}}{\varepsilon^n}=0\,.~\Box$$

\medskip
Ещё одно утверждение может быть получено в случае, если $Q\in
L_{loc}^{\alpha}({\Bbb B}^n).$

\medskip
\begin{theorem}\label{th2}{\sl\,
Пусть $n\geqslant 2,$ $n<p<\infty,$ $f:{\Bbb B}^n\rightarrow {\Bbb
R}^n$ -- открытое дискретное кольцевое $Q$-отоб\-ра\-же\-ние
относительно $p$-модуля в нуле, $f(0)=0.$ Предположим, что $Q:{\Bbb
B}^n\rightarrow[0,\,\infty]$ -- локально интегрируемая функция в
${\Bbb B}^n$ в степени $\alpha> 1.$ Пусть $K\subset {\Bbb B}^n$ --
произвольный компакт, удовлетворяющий условию $0\in {\rm Int}\, K.$
Тогда имеет место оценка:
$$
\limsup\limits_{x\rightarrow
0}\frac{|f(x)|}{|x|^{1+\frac{n}{\alpha(p-n)}}}\geqslant C>0\,,
$$
где $C$  -- некоторая положительная постоянная, зависящая только от
размерности пространства $n,$ $p,$ $\alpha$ и компакта $K.$}
\end{theorem}

\begin{proof}
Выберем произвольным образом компакт $K\subset {\Bbb B}^n,$
удовлетворяющий условию $0\in {\rm Int}\, K.$ Поскольку $Q\in
L_{loc}^{\alpha}({\Bbb B}^n),$ найдётся постоянная
$\overline{C}=\overline{C}(K)<\infty$ такая, что
$\int\limits_KQ^{\alpha}(x)\,dm(x)\leqslant \overline{C}(K).$
Повторяя теперь рассуждения, проведённые при доказательстве теоремы
\ref{th1}, мы вновь получаем соотношения (\ref{eq101}) и
(\ref{eq3}). Кроме того, поскольку по выбору компакта $K$ точка 0
является его внутренней точкой, при достаточно малых $\varepsilon>0$
кольцо $A(0, \varepsilon, 2\varepsilon)$ лежит в $K.$

\medskip Оценим теперь интеграл справа в (\ref{eq101}) сверху,
используя неравенство Гёльдера с показателями $\alpha$ и
$\alpha^{\,\prime}=\frac{\alpha}{\alpha-1}>1,$
$1/{\alpha}+1/{\alpha^{\,\prime}}=1.$ Учитывая сказанное выше, будем
иметь:
$$\int\limits_{A(0, \varepsilon,
2\varepsilon)}Q(x)\,dm(x)\leqslant
\left(\int\limits_KQ^{\,\alpha}(x)dm(x)\right)^{1/\alpha}\cdot
(2\Omega^{1/n}_n\varepsilon)^{\frac{n(\alpha-1)}{\alpha}}\,\leqslant
C_1(K)\cdot \varepsilon^{\frac{n\alpha-n}{\alpha}}\,,
$$
где $C_1=C_1(K)$ -- некоторая новая постоянная, зависящая только от
функции $Q,$ компакта $K,$ $n$ и степени $\alpha.$ Применяя
(\ref{eq101}) и (\ref{eq3}), мы будем иметь:
\begin{equation}\label{eq4} (m(f(B(0, 2\varepsilon)))^{1/n}
\geqslant C_2\cdot
\varepsilon^{(\frac{n\alpha-n}{\alpha}-p)\cdot\frac{1}{n-p}}=C_2\cdot
\varepsilon^{\frac{n\alpha-n-\alpha p}{\alpha n-\alpha p}}=C_2\cdot
\varepsilon^{1+\frac{n}{\alpha(p-n)}}\,,\end{equation}
где $C_2$ -- некоторая положительная постоянная, зависящая только от
функции $Q,$ компакта $K,$ $n$ и степени $\alpha.$
Из (\ref{eq4}) вытекает, что
\begin{equation}\label{eq5} \frac{(m(f(B(0, 2\varepsilon)))^{1/n}}{\varepsilon^{1+\frac{n}{\alpha(p-n)}}}
\geqslant C_2>0\,.\end{equation}
Полагая $L_f(\varepsilon)=\max\limits_{|x|=\varepsilon}|f(x)|$ и
используя соотношения (\ref{eqA28}) и (\ref{eq5}), мы получим, что
$$\limsup\limits_{x\rightarrow
0}\frac{|f(x)|}{|x|^{1+\frac{n}{\alpha(p-n)}}}=\limsup\limits_{\varepsilon\rightarrow
0}\frac{L_f(2\varepsilon)}{(2\varepsilon)^{1+\frac{n}{\alpha(p-n)}}}
\geqslant\limsup\limits_{\varepsilon\rightarrow 0}\frac{(m(f(B(0,
2\varepsilon))))^{1/n}}{(2\varepsilon)^{1+\frac{n}{\alpha(p-n)}}}\geqslant$$
%
%\begin{equation}\label{eq6}
$$\geqslant C_2/2^{1+\frac{n}{\alpha(p-n)}}\,.$$
%\end{equation}
%
Полагая $C:=C_2/2^{1+\frac{n}{\alpha(p-n)}},$ получаем необходимое
заключение.~$\Box$
\end{proof}

\medskip
Заметим, что от условия $Q\in L^{\alpha}_{loc}({\Bbb B}^n),$
участвующего в теореме \ref{th2}, вообще говоря, нельзя отказаться,
что составляет содержательную часть следующего утверждения.

\medskip
\begin{theorem}\label{th3}
{\sl Для произвольного $\alpha>1$ найдётся кольцевой
$Q$-гомеоморфизм $f:{\Bbb B}^n\rightarrow {\Bbb B}^n$ в точке
$x_0=0$ относительно $p$-модуля, $p>n,$ такой что $Q\not\in
L^{\alpha}_{loc}({\Bbb B}^n)$ и, при этом,
\begin{equation}\label{eq12A} \lim\limits_{x\rightarrow
0}\frac{|f(x)|}{|x|^{1+\frac{n}{\alpha(p-n)}}}=0\,. \end{equation}}
\end{theorem}

\medskip
\begin{proof}Зафиксируем произвольным образом $\alpha>1,$ $\varepsilon>0$ и
$p>n.$ Положим
\begin{equation}\label{eq10A}
f(x)=x\cdot|x|^{\frac{n}{\alpha(p-n)}+\varepsilon}\,.
\end{equation}
Из определения вытекает, что $f$ -- гомеоморфизм единичного круга в
себя, $f(0)=0,$ кроме того,
$f^{-1}(y)=\frac{y}{|y|}\cdot|y|^{1/\beta},$ где
$\beta=\frac{n}{\alpha(p-n)}+\varepsilon+1.$ Отсюда следует, что $f,
f^{-1}\in C^1({\Bbb B}^n\setminus\{0\});$ в частности, $f$
дифференцируемо почти всюду, обладает $N$ и $N^{\,-1}$-свойствами
Лузина и, кроме того, $f^{\,-1}$ абсолютно непрерывно на $p$-почти
всех кривых. Тогда по \cite[теорема~1.1]{SalSev$_3$} отображение $f$
является кольцевым $Q$-отображением в каждой точке $x_0\in{\Bbb
B}^n$ при
$$Q(x):=K_{I, p}(x,f)\quad =\quad\left\{
\begin{array}{rr}
\frac{J(x,f)}{l^p(x, f)}, & J(x,f)\ne 0,\\
1,  &  f^{\,\prime}(x)=0, \\
\infty, & \text{в\,\,остальных\,\,случаях},
\end{array}
\right.$$
где $l\left(f^{\,\prime}(x)\right)\,=\,\,\,\min\limits_{h\in {\Bbb
R}^n \setminus \{0\}} \frac {|f^{\,\prime}(x)h|}{|h|}$ и $J(x,
f)={\rm det}\, f^{\,\prime}(x).$ Подсчитаем $K_{I, p}(x,f),$ для
чего воспользуемся предложением 5.1 в \cite{IS}. В обозначениях
этого предложения
$$\lambda_r(x)=\frac{\partial |f(x)|}{\partial|x|}=\left(\frac{n}{\alpha(p-n)}+\varepsilon+1\right)|x|^{\frac{n}{\alpha(p-n)}+\varepsilon}\,,$$
$$\lambda_{\tau}(x)=\frac{|f(x)|}{|x|}=|x|^{\frac{n}{\alpha(p-n)}+\varepsilon}\,.$$
Заметим, что $l\left(f^{\,\prime}(x)\right)=\min\{\lambda_r(x),
\lambda_{\tau}(x)\}$ и $|J(x,
f)|=\lambda^{n-1}_{\tau}(x)\cdot\lambda_r(x)$ (см.
\cite[разд.~5.1]{IS}). Очевидно, $\lambda_{\tau}(x)<\lambda_r(x),$
поэтому
$$K_{I, p}(x,f)=\left(\frac{n}{\alpha(p-n)}+\varepsilon+1\right)\cdot
\frac{|x|^{(n-1)(\frac{n}{\alpha(p-n)}+\varepsilon)}\cdot
|x|^{\frac{n}{\alpha(p-n)}+\varepsilon}}
{|x|^{p(\frac{n}{\alpha(p-n)}+\varepsilon)}}\,=$$
$$=\left(\frac{n}{\alpha(p-n)}+\varepsilon+1\right)\cdot|x|^{-n/\alpha+\varepsilon(n-p)}\,=C\cdot|x|^{-n/\alpha+\varepsilon(n-p)}\,,$$
где $C=\frac{n}{\alpha(p-n)}+\varepsilon+1.$ Пусть $K$ --
произвольный компакт в ${\Bbb B}^n$ такой, что $B(0,
\varepsilon_0)\subset K$ при некотором $0<\varepsilon_0<1.$ Тогда по
теореме Фубини
$$\int\limits_{K}K^{\alpha}_{I, p}(x,f)dm(x)\geqslant \int\limits_{B(0, \varepsilon_0)}K^{\alpha}_{I, p}(x,f)\,dm(x)=$$
\begin{equation}\label{eq11}=C^{\alpha}\cdot\int\limits_0^{\varepsilon_0}\int\limits_{S(0,
r)}|x|^{-n+\alpha\varepsilon(n-p)}dSdr=C^{\alpha}\omega_{n-1}
\int\limits_0^{\varepsilon_0}r^{n-1}r^{-n+\alpha\varepsilon(n-p)}dr=
\end{equation}
$$=C^{\alpha}\omega_{n-1}
\int\limits_0^{\varepsilon_0}r^{\alpha\varepsilon(n-p)-1}dr=\infty\,,$$
так как показатель степени при $r$ меньше $-1.$ Очевидно также, что
соотношение (\ref{eq12A}) имеет место, т.е., заключение теоремы
\ref{th2} не выполнено. Причиной последнего является расходимость
интеграла слева в (\ref{eq11}).~$\Box$
\end{proof}

\medskip
КОНТАКТНАЯ ИНФОРМАЦИЯ

\medskip
\noindent{{\bf Руслан Радикович Салимов} \\
Институт математики НАН Украины \\
ул. Терещенковская, д. 3 \\
г. Киев-4, Украина, 01 601\\
тел. +38 095 630 85 92 (моб.), e-mail: ruslan.salimov1@gmail.com}

\medskip
\noindent{{\bf Евгений Александрович Севостьянов, Антонина Александровна Маркиш} \\
Житомирский государственный университет им.\ И.~Франко\\
ул. Большая Бердичевская, 40 \\
г.~Житомир, Украина, 10 008 \\ тел. +38 066 959 50 34 (моб.),
e-mail: esevostyanov2009@gmail.com, tonya@bible.com.ua}

\end{document}